\documentclass[12pt,leqno]{article}

\usepackage{amsfonts}
\usepackage{amssymb}
\usepackage[all]{xy}

\newcommand{\cM}{{\cal M}}

\newcommand{\cO}{{\cal O}}

\newcommand{\Aut}{{\rm Aut}}

\newcommand{\Spec}{{\rm Spec}}
\newcommand{\Proj}{{\rm Proj}}

\newcommand{\Crit}{{\rm Crit}}

\newcommand{\CC}{{\mathbb C}}

\newcommand{\NN}{{\mathbb N}}

\newcommand{\PP}{{\mathbb P}}
\newcommand{\QQ}{{\mathbb Q}}

\newcommand{\ra}{\rightarrow}
\def\rightepi{{\longrightarrow \kern-0.7em \rightarrow}}

\newcommand{\notteilt}{{\,\not{\kern-0.075em|}\,}}
\def\antiddots{\mathinner{\mkern1mu\raise1pt\vbox{\kern7pt\hbox{.}}\mkern2mu
    \raise4pt\hbox{.}\mkern2mu\raise7pt\hbox{.}\mkern1mu}}

\begin{document}

\vspace*{15ex}

\begin{center}
{\LARGE\bf Belyi's Theorem Revisited}\\
\bigskip
by\\
\bigskip
{\sc Bernhard K\"ock}
\end{center}

\bigskip

\begin{quote}
{\footnotesize {\bf Abstract}. We give an elementary,
self-contained and quick proof of Belyi's theorem. As a by-product
of our proof we obtain an explicit bound for the degree of the
defining number field of a Belyi surface.


}

\end{quote}

\bigskip

\section*{Introduction}

The main purpose of this paper is to give an elementary,
self-contained and quick proof of the following famous theorem by
Belyi (see \cite{Bel}).

{\bf Theorem.} {\em A complex smooth projective curve $X$ is
defined over a number field, if and only if there exists a
non-constant morphism $t: X \ra \PP^1_\CC$ with at most 3 critical
values.}

While the only-if-direction is just a fairly elementary and short
algorithm which is well explained in the literature and, once
more, in Lemmas (3.4) through (3.6) below, I found it difficult to
understand the proofs of the if-direction existing in the
literature. So, the main focus in this paper is on the
if-direction which is also called the ``obvious part" which
somebody familiar with the results of Weil's paper \cite{We}, in
particular Theorem 4, and with the mathematical language used
there may consider as justified.

As already observed by Wolfart in his paper \cite{Wo}, the notion
{\em moduli field} allows an elegant way to split up the
if-direction into two assertions. However, rather than using the
(absolute) moduli field of a complex smooth projective curve $X$
as in \cite{Wo}, we will use the (relative) moduli field of a
finite morphism $t: X \ra \PP^1_\CC$ which, by definition, is the
subfield $\CC^{U(X,t)}$ of $\CC$ fixed by the subgroup $U(X,t)$ of
all automorphisms $\sigma$ of $\CC$ such that there is an
isomorphism between the curve $X^\sigma$ and $X$ compatible with
the covering $t$, see Notation~(1.1) and Definition~(2.1). We will
prove the following two assertions which obviously imply the
if-direction in Belyi's theorem.

Let $X$ be a complex smooth projective curve, and let
$t: X \ra \PP^1_\CC$ be a finite morphism. Then we have:\\
\hspace*{2ex}(a) If the critical values lie in $\{0,1,\infty\}$,
then the moduli field of $t$ is a number field.\\ \hspace*{2ex}(b)
$X$ and $t$ are defined over a finite extension of the moduli
field of $t$.

Assertion~(a), a special case of Corollary~(3.2), follows from the
fact (see Proposition~(3.1)) that there are at most finitely many
isomorphism classes of coverings $t:X \ra \PP^1_\CC$ of given
degree $d$ and given subset $S$ of $\PP^1_\CC$ of critical values.
This fact occurs implicitly at several places in the literature.
We include a short, self-contained and elementary proof which, in
contrast to the existing literature, avoids any non-standard or
highly sophisticated notion or fact. A slightly strengthened and
generalized version of Assertion~(b) will be given in
Theorem~(2.2). Its proof is based on ideas going back to
Grothendieck and Coombes/Harbater. Apart from being the most
original part in our proof of Belyi's theorem, it also yields an
interesting explicit bound for the degree of the defining number
field which seems to be new, see Corollary~(3.7).

As indicated already above, it is obviously possible to replace,
in the proof of the if-direction, Assertion~(b) by the following
absolute analogue: Any complex smooth projective curve is defined
over a finite extension of its moduli field. This assertion is of
independent interest and can be strengthened in many cases, see
Example~(1.7) and Corollary~(1.11). However, it is probably fair
to say that its proof is much more sophisticated, see Wolfart's
paper \cite{Wo} or the more recent paper \cite{HH} by Hammer and
Herrlich. (In contrast to \cite{Wo}, the proof in the latter paper
does not use Weil's language of generic points and also works for
ground fields of positive characteristic.)

{\bf Acknowledgments}. I would like to thank my colleagues at
Karlsruhe and Southampton for their encouraging interest and for
helpful hints. Furthermore, I would like to thank J.\ Wolfart for
extensive e-mail discussions. Finally, I would like to thank the
referees for some helpful comments improving the presentation of
the paper.

\bigskip

\section*{1. The Moduli Field of a Curve}

In this paper, a {\em curve over a field $C$} means a smooth
projective geometrically connected variety of dimension $1$ over
$C$, and a {\em variety over $C$} is an integral separated scheme
$X$ together with a morphism $p: X \ra \Spec(C)$ of finite type.
(For the purposes of this paper it is convenient and appropriate
to describe everything in the language of schemes, but we will not
need any deeper insight into the theory of schemes.  In
particular, a reader familiar only with the language of classical
varieties will presumably be able to read this paper without any
difficulties.)

{\bf (1.1) Notation.} Let $C$ be a field and let $p: X \ra
\Spec(C)$ be a variety. For any $\sigma \in \Aut(C)$, we denote by
$X^\sigma/C$ the variety consisting of the scheme $X$ and the
structure morphism $X \,\, \stackrel{p}{\longrightarrow} \,\,
\Spec(C) \,\, \stackrel{\Spec(\sigma)}{\longrightarrow} \,\,
\Spec(C)$.

Note that the scheme underlying the variety $X^\sigma/C$ is the
same as the scheme underlying the variety $X/C$. In particular,
the function field of $X^\sigma/C$ is the same as the function
field of $X/C$, and an isomorphism between $X^\sigma/C$ and $X/C$
is an {\em auto}morphism of the scheme $X$. Unfortunately, the
concept of changing the structure morphisms by $\Spec(\sigma)$
does not exist in the language of classical varieties. However,
the following remark shows that the variety $X^\sigma/C$ is
isomorphic to that variety which is usually denoted by
$X^\sigma/C$ in the language of classical varieties; the
isomorphism is induced by $\sigma$.

{\bf (1.2) Remark.} Let $C$ be a field, let $\sigma \in \Aut(C)$
and let $X$ be a subvariety of $\PP^n_C$ given by the homogeneous
polynomials $f_1, \ldots, f_m \in C[X_0, \ldots, X_n]$, i.e., $X =
V(f_1, \ldots, f_m)$. Let $\sigma$ also denote the induced
automorphism of $C[X_0, \ldots, X_n]$. Then the variety
$X^\sigma/C$ is given by the polynomials $\sigma^{-1}(f_1),
\ldots, \sigma^{-1}(f_m) \in C[X_0, \ldots, X_n]$.

{\em Proof.} Let $\bar{\sigma}$ denote the isomorphism between
$C[X_0, \ldots, X_n]/(\sigma^{-1}(f_1), \ldots, \sigma^{-1}(f_m))$
and $C[X_0, \ldots, X_n]/(f_1, \ldots f_m)$ induced by $\sigma$.
Then ${\rm Proj}(\bar{\sigma})$ is the desired isomorphism between
the varieties $V(f_1, \ldots, f_m)^\sigma/C$
%\longrightarrow \Spec(C)
%\,\, \stackrel{\Spec(\sigma)}{\longrightarrow} \,\, \Spec(C)$
and $V(\sigma^{-1}(f_1), \ldots, \sigma^{-1}(f_m))/C$.
\hspace*{\fill} $\Box$
%\longrightarrow \Spec(C)$.

For readers not familiar with the language of schemes, the
following explanation might also help to understand the notation
$X^\sigma/C$. It is well-known and we will frequently use that a
curve $X/C$ is the same as a finitely generated field $K$ of
transcendence degree 1 over $C$ such that $C$ is algebraically
closed in $K$. It is important here that the embedding of $C$ into
$K$ belongs to the notion of a curve. Changing this embedding by
an automorphism $\sigma$ of $C$ yields a new curve which
corresponds to the notation $X^\sigma/C$ in (1.1).

From now on we assume that $C$ is an algebraically closed field of
characteristic $0$.

{\bf (1.3) Definition.} The {\em moduli field} of a variety $X/C$
is the field $M(X) := C^{U(X)}$ fixed by the subgroup
\[ U(X) := \{\sigma \in \Aut(C): X^\sigma/C \textrm{ is isomorphic
to }X/C\}\] of $\Aut(C)$.

As usual, we say that a variety $X/C$ {\em is defined over the
subfield $K$ of $C$}, iff there is a variety $X_K/K$ such that
$X/C$ is isomorphic to $X_K \times_K C/C$, i.e., iff $X/C$ can be
covered by affine varieties which are given by polynomials with
coefficients in $K$. In this case, the subgroup $\Aut(C/K)$ of
$\Aut(C)$ is obviously contained in $U(X)$, hence the moduli field
$M(X)$ is contained in $K$ by the following folklore lemma (which
at the same time is a central argument in the proof of Belyi's
theorem). In particular, if $X/C$ is defined over its moduli field
$M(X)$, then $M(X)$ is the smallest field of definition for $X/C$.

{\bf (1.4) Lemma.} {\em Let $K$ be a subfield of $C$. Then, any
automorphism of $K$ can be extended to an automorphism of $C$.
Furthermore, we have:
\[ C^{\Aut(C/K)} = K.\]}

{\em Proof.} The first assertion is well-known and easy to prove.
The inclusion $K \subseteq C^{\Aut(C/K)}$ is a tautology. The
reverse inclusion is equivalent to the assertion that, for any $x
\in C\backslash K$, there is a $\sigma \in \Aut(C/K)$ with
$\sigma(x) \not= x$. If $x$ is transcendent over $K$, then mapping
$x$, for instance, to $-x$ yields a $K$-automorphism of $K(x)$
which does not fix $x$. This automorphism can be extended to the
desired automorphism $\sigma$ of $C$ by the first assertion. If
$x$ is algebraic over $K$, we choose a $y \in C\backslash \{x\}$
which is $K$-conjugate to $x$. Then, mapping $x$ to $y$ yields a
$K$-embedding of $K(x)$ into the normal closure~$L$ of $K(x)$ over
$K$. This embedding can be extended to a $K$-automorphism of $L$
and then, again by the first assertion, to the desired
$K$-automorphism $\sigma$ of $C$. \hspace*{\fill} $\Box$

We call a subgroup $U$ of $\Aut(C)$ {\em closed}, iff there is a
subfield $K$ of $C$ with $U = \Aut(C/K)$. Lemma~(1.4) implies that
we have a bijective Galois correspondence between the set of
subfields of $C$ and the set of closed subgroups of $\Aut(C)$. In
particular we have $U= \Aut(C/C^U)$ for any closed subgroup $U$ of
$\Aut(C)$. For any field $C$ as above, there exist non-closed
subgroups of $\Aut(C)$ (even of finite index); in the case
$C=\bar{\QQ}$ this is a well-known fact in infinite Galois theory;
in the general case, the preimage of a non-closed subgroup of
$\Aut(\bar{\QQ})$ under the canonical epimorphism $\Aut(C) \ra
\Aut(\bar{\QQ})$ is a non-closed subgroup of $\Aut(C)$. It follows
from Lemma~(1.5) and Theorem (1.8) below that the subgroup $U(X)$
of $\Aut(C)$ introduced in Definition (1.3) is closed, if $X/C$ is
a curve.

{\bf (1.5) Lemma.} {\em Let $U$ be a subgroup of $\Aut(C)$ such
that there is a finite field extension $K/C^U$ with $\Aut(C/K)
\subseteq U$. Then $U$ is closed.}

{\em Proof.} We may assume that $K/C^U$ is a finite Galois
extension. Then $C^U$ is the field fixed by the image $B$ of
$U/\Aut(C/K)$ under the canonical isomorphism
\[\Aut(C/C^U)/\Aut(C/K) \,\, \tilde{\ra} \,\, \Aut(K/C^U).\] Thus,
$B= \Aut(K/C^U)$, and hence $U = \Aut(C/C^U)$ is closed.
\hspace*{\fill} $\Box$

For later purposes we record the following lemma.

{\bf (1.6) Lemma.} {\em Let $U$ be a subgroup of $\Aut(C)$ and let
$V$ be a subgroup of $U$ of finite index. Then the field extension
$C^V/C^U$ is finite. If $V$ is a normal subgroup of $U$ or if $U$
is closed, then we have $[C^V: C^U] \le [U:V]$. If $V$ is closed,
then we even have $[C^V:C^U] = [U:V]$.}

{\em Proof.} It is easy to see and well-known that there is a
normal subgroup $W$ of $U$ of finite index which is contained in
$V$. Then we obviously have a canonical homomorphism $U/W \ra
\Aut(C^W/C^U)$. The field fixed by the image $B$ of this
homomorphism is $C^U$. Thus $C^W/C^U$ is a finite Galois extension
and we have $B = \Aut(C^W/C^U)$. Hence we obtain:
\[[C^W:C^U] = {\rm ord}(\Aut(C^W/C^U)) \le {\rm ord}(U/W) =
[U:W].\] This implies the first assertion of Lemma~(1.6) and also
the second assertion in the case that already $V$ is a normal
subgroup of $U$. Furthermore we have:
\begin{eqnarray*}
\lefteqn{[C^V:C^U] = \frac{[C^W:C^U]}{[C^W:C^V]} = \frac{{\rm
ord}(\Aut(C^W/C^U))}{{\rm ord}(\Aut(C^W/C^V))} =}\\&& =
\left|\frac{\Aut(C^W/C^U)}{\Aut(C^W/C^V)}\right| =
\left|\frac{\frac{\Aut(C/C^U)}{\Aut(C/C^W)}}
{\frac{\Aut(C/C^V)}{\Aut(C/C^W)}} \right| =
\left|\frac{\Aut(C/C^U)}{\Aut(C/C^V)}\right|.
\end{eqnarray*}
This implies the second assertion in the case that $U$ is closed
and also the third assertion, because, if $V$ is closed, then also
$U$ is closed by Lemma~(1.5) and by the first assertion. Thus,
Lemma~(1.6) is proved. \hspace*{\fill} $\Box$

The rest of this section deals with the question how far the
moduli field of a curve $X/C$ is away from being a field of
definition for $X/C$. It is not used in the proof of Belyi's
theorem. We start with the following well-known elementary
example.

{\bf (1.7) Example.} {\em Let $X/C$ be a curve of genus $0$ or
$1$. Then $X/C$ is defined over its moduli field $M(X)$.}

{\em Proof.} In the case $g=0$, $X$ is isomorphic to the
projective line which is defined over $\QQ$, the smallest field of
characteristic $0$. Now let $g=1$. Then $X/C$ is an elliptic
curve. Let $j\in C$ denote the $j$-invariant of $X/C$. Then we
have $U(X) = \{\sigma \in \Aut(C): \sigma(j) =j\} =
\Aut(C/\QQ(j))$, thus $M(X) = \QQ(j)$ by Lemma~(1.4). Furthermore
it is well-known that $X/C$ is defined over $\QQ(j)$ (e.g., see
Proposition~1.4 on p.~50 in \cite{Si}). \hspace*{\fill} $\Box$

In general we have:

{\bf (1.8) Theorem.} {\em Let $X/C$ be a curve. Then $X/C$ is
defined over a finite extension of its moduli field $M(X)$.}

{\em Proof.} See Theorem 4 in \cite{Wo} or Theorem~5 in \cite{HH}.
\hspace*{\fill} $\Box$

The object of the following considerations is to strengthen
Theorem~(1.8). The basic tool for this is the following theorem
which is a slight weakening of Theorem 1 in Weil's paper
\cite{We}. I hope that the given formulation and the given proof
make this theorem easier to access.

{\bf (1.9) Theorem.} {\em Let $L$ be a field, let $G$ be a finite
subgroup of $\Aut(L)$ and let $X/L$ be a variety. We suppose that,
for any $\sigma \in G$, we are given a birational map $f_\sigma:
X^\sigma \ra X$ over $L$ such that
\[f_{\sigma \tau} = f_\sigma \circ f_\tau^\sigma \textrm{ for all }
\sigma, \tau \in G.\] Then there is a variety $X_K$ over the fixed
field $K:= L^G$ such that $X_K \times_K L/L$ is birationally
equivalent to $X/L$.}

Here, the notation $f_\tau^\sigma$ means that we consider the
(auto)morphism $f_\tau$ (defined on some open subscheme of $X$) as
a rational morphism from the variety $X^{\sigma \tau} =
(X^\tau)^\sigma$ to the variety $X^\sigma$. In the language of
classical varieties, the notation $f_\tau^\sigma$ means that we
apply $\sigma^{-1}$ to the polynomials defining $f_\tau$ (cf.\
Remark (1.2)).

{\em Proof.} We have to show that there is a finitely generated
field $V$ over $K$ such that $L \otimes_K V$ is $L$-isomorphic to
the function field $W$ of $X$. This follows from the following
lemma applied to the action $G \ra \Aut(W)$, $\sigma \mapsto
f_\sigma^*$, of $G$ on $W$. Note that the fixed field $V := W^G$
is finitely generated over $K$ (being an intermediate field of
$W/K$). \hspace*{\fill} $\Box$

{\bf (1.10) Lemma} (Galois descent). {\em Let $L$ be a field, let
$G$ be a finite subgroup of $\Aut(L)$, and let $W$ be a vector
space over $L$ together with a semilinear action of $G$ on $W$
(i.e., $\sigma(aw) = \sigma(a) \sigma(w)$ for all $\sigma \in G$,
$a\in L$, and $w\in W$). We set $K:=L^G$. Then the following
canonical $L$-homomorphism is bijective:
\[L \otimes_K W^G \,\, \tilde{\longrightarrow}\,\, W.\]}

{\em Proof.} The proof of the injectivity is rather
straightforward and the surjectivity follows from the linear
independence of characters (see books on Galois theory for
details). \hspace*{\fill} $\Box$

We recall that the automorphism group $\Aut(X/C)$ of a curve $X/C$
of genus $g \ge 2$ is finite (see Exercise 5.2 on p.\ 348 in
\cite{Ha}). If $g\ge 3$, then $\Aut(X/C)$ is even ``generically
trivial" (see Exercise 5.7 on p.\ 348 in \cite{Ha}); in
particular, the following corollary implies that ``almost all"
curves $X/C$ of genus $g\ge 3$ are defined over their moduli
field. If $C=\bar{\QQ}$, this corollary is a special case of
Theorem 3.1 in the paper \cite{DE} by D\`ebes and Emsalem. If
$C=\CC$, it is mentioned in Wolfart's paper \cite{Wo}, but without
proof. We here give a complete proof.

{\bf (1.11) Corollary.} {\em Let $X/C$ be a curve of genus $g \ge
2$. Then the quotient curve $X/\Aut(X/C)$ is defined over $M(X)$.}

{\em Proof.} By Theorem (1.8), there is a model $X_{L'}/L'$ of
$X/C$ over a finite Galois extension $L'$ of $M:=M(X)$. By Lemma
(1.5), we have $U(X) = \Aut(C/M)$ and the canonical homomorphism
$U(X) \ra \Aut(L'/M)$ is surjective. For any $\tau \in
\Aut(L'/M)$, we choose a preimage $\tilde{\tau} \in U(X)$ and an
isomorphism $f_{\tilde{\tau}}: X^{\tilde{\tau}} \ra X$ of curves
over $C$. By Lemma~(1.12) below, there is a Galois extension $L$
of $M$ such that $L'$ is contained in $L$ and such that all
isomorphisms $f_{\tilde{\tau}}$, $\tau \in \Aut(L'/M)$, and  all
automorphisms of $X/C$ are defined over $L$. We set $X_L := X_{L'}
\times_{L'} L$, $G:= \Aut(L/M)$, and write $\bar{\sigma}$ for the
image of $\sigma \in G$ in $\Aut(L'/M)$ and $g_{{\sigma}}$ for the
isomorphism $X_L^\sigma \,\, \tilde{\ra}\,\, X_L$ with
$g_{{\sigma}} \times_L C = f_{\tilde{\bar{\sigma}}}$. The
isomorphism $g_{{\sigma}}$ induces an isomorphism
\[h_\sigma: X_L^\sigma/\Aut(X_L^\sigma/L) \ra X_L/\Aut(X_L/L)\]
between the quotient curves which does not depend on the choice of
the isomorphism $f_{\tilde{\bar{\sigma}}}:
X^{\tilde{\bar{\sigma}}} \ra X$. In particular, we have $h_{\tau
\sigma} = h_{\sigma} \circ h_\tau$ for all $\sigma, \tau \in G$.
By Theorem 1.9, the curve $X_L/\Aut(X_L/L)$ and hence the curve
$X/\Aut(X/C)$ is defined over $M=L^G$. \hspace*{\fill} $\Box$

I would like to thank J.\ Wolfart for the elementary main idea in
the proof of the following folklore fact from Algebraic Geometry.

{\bf (1.12) Lemma.} {\em Let $N$ be an algebraically closed
subfield of $C$ and let $X/N$ and $Y/N$ be curves of genus $\ge
2$. Then, any isomorphism between $X_C:= X \times_N C$ and $Y_C:=
Y \times_N C$ is already defined over $N$. In particular, the
following canonical homomorphism is bijective:
\[ \Aut(X/N) \,\, \tilde{\ra} \,\, \Aut(X_C/C).\]}

{\em Proof.} We choose $t_1, t_2 \in K(Y)\backslash N$ with $K(Y)
= N(t_1, t_2)$ and denote the minimal polynomial of $t_2$ over
$N(t_1)$ by $g \in N(t_1)[T_2]$. By the usual dictionary between
curves and function fields, we then have a natural bijection
between the set of isomorphisms from $X_C/C$ to $Y_C/C$ and the
set $\cM$ of pairs $(\alpha_1, \alpha_2)$ in $K(X_C) \backslash C$
with $K(X_C) = C(\alpha_1, \alpha_2)$ and $g(\alpha_1, \alpha_2) =
0$. Since the genus of $X$ and $Y$ is greater than or equal to
$2$, the set $\cM$ is finite. On the other hand, if the set $\cM$
contains $(\alpha_1,\alpha_2)$, then it also contains
$(\tau(\alpha_1), \tau(\alpha_2))$ for any $\tau \in
\Aut(K(X_C)/N)$ with $\tau(C) = C$. Furthermore, the set
$\{\sigma(x) : \sigma \in \Aut(C/N)\}$ is infinite, if $x \in
C\backslash N$ (see the proof of Lemma~(1.4)), and any $\sigma \in
\Aut(C/N)$ can obviously be extended to a $\tau \in
\Aut(K(X_C)/N)$. Thus, any pair $(\alpha_1, \alpha_2)$ as above
comes already from $K(X) \backslash N$. So, any isomorphism
between $X_C$ and $Y_C$ is already defined over $N$.
\hspace*{\fill} $\Box$

\bigskip

\section*{2. The Moduli Field of a Covering}

Let $C$ be an algebraically closed field of characteristic $0$ and
let $t:X \ra \PP^1_C$ be a finite morphism from a curve $X/C$ to
the projective line $\PP^1_C$. We will denote the degree of $t$ by
$\deg(t)$ and we will use term {\em critical value} for any point
$Q \in \PP^1_\CC$ which has less than $\deg(t)$ preimages under
$t$.

{\bf (2.1) Definition.} The {\em moduli field of $t$} is the field
$M(X,t):= C^{U(X,t)}$ fixed by the subgroup $U(X,t)$ of $U(X)$
consisting of all $\sigma \in \Aut(C)$ such that there exists an
isomorphism $f_\sigma: X^\sigma \ra X$ of varieties over $C$ such
that the following diagram commutes:
\[ \xymatrix{
X^\sigma \ar[d]^{t^\sigma} \ar[rr]^{f_\sigma} && X \ar[d]^{t} \\
(\PP^1_C)^\sigma \ar[rr]^{\Proj(\sigma)} && \PP^1_C;} \] here,
$\Proj(\sigma)$ means the automorphism of the scheme $\PP^1_C =
\Proj(C[T_0,T_1])$ induced by the extension of the automorphism
$\sigma \in \Aut(C)$ to $C[T_0,T_1]$ (denoted $\sigma$ again).

Obviously we have $M(X) \subseteq M(X,t)$. The following theorem
is the analogue of Theorem (1.8) but much easier to prove. We will
merely use the Riemann-Roch theorem for curves and basic facts of
the ramification theory for curves. More precisely, we combine
some ideas of the proof of Proposition 2.1 on p.~10 in the book
\cite{MM} by Malle and Matzat (going back to Grothendieck) with an
idea by Coombes and Harbater (see Proposition 2.5 on p.~830 in
\cite{CH}). In fact, if $C = \CC$, the second assertion of
Theorem~(2.2) is Proposition 2.5 in \cite{CH}.

{\bf (2.2) Theorem.} {\em The curve $X/C$ and the morphism $t$ are
defined over a finite extension of $M(X,t)$. If $t$ is a Galois
covering (i.e., if the corresponding extension of function fields
is Galois), then $X/C$ and $t$ are defined over $M(X,t)$ itself.}

{\em Proof.} We choose a $\QQ$-rational point $Q$ of $\PP^1_C$
which is not a critical value of $t$, and we choose a point $P$ in
the fibre $t^{-1}(Q)$. By the theorem of Riemann-Roch (see
Theorem~1.6 on p.\ 362 in \cite{Ha}) applied to the divisor $D:=
({\rm genus}(X)+1)[P]$, there is a meromorphic function $z \in
K(X)\backslash C$ such that $P$ is the only pole of $z$. Then we
have $K(X) = C(t,z)$ where, here, $t$ is considered as a
meromorphic function on $X$; for the field extension
$K(X)/C(t,z)$ is a subextension of $K(X)/C(t)$ and of $K(X)/C(z)$,
hence the corresponding morphism of curves is both unramified and
totally ramified at $P$. We assume furthermore that we have chosen
$z$ in such a way that the pole order $m:= -{\rm ord}_P(z) \in
\NN$ is minimal. Then we have
\[V:= \{x \in K(X): {\rm ord}_P(x) \ge -m\ \textrm{ and }
{\rm ord}_Q(x) \ge 0 \textrm{ for all } Q \in X \backslash \{P\}\} = C \oplus Cz;\]
for, for any $x_1, x_2 \in V$ with ${\rm ord}_P(x_i) = -m$,
$i=1,2$, there is a constant $\alpha \in C$ with $-{\rm ord}_P(x_1
- \alpha x_2) < m$, and then $x_1 - \alpha x_2$ is a constant
function, since $m$ was minimal. By the choice of $Q$, the
meromorphic function $t-Q$ on $X$ is a local parameter on $X$ in
$P$; if $C=\CC$, this means, in the language of Riemann surfaces,
that $t-Q$ yields a chart of $X(\CC)$ in a neighborhood of $P$
which maps $P$ to $0$. There is obviously a unique function $z'
\in V$ such that the leading coefficient (i.e., the coefficient of
$(t-Q)^{-m}$) and the constant coefficient (i.e., the coefficient
of $(t-Q)^0$) in the Laurent expansion of $z'$ with respect to the
local parameter $t-Q$ are equal to $1$ and $0$, respectively. (In
the language of Algebraic Geometry, the term ``Laurent expansion
of $z'$" means ``the image of $z'$ in the quotient field of the
completion $\hat{\cO}_{X,P} = C[[t-Q]]$ of the local ring
$\cO_{X,P}$"). We may and we will assume that $z =z'$. We now
claim that the minimal polynomial of $z$ over $C(t)$ has
coefficients in $k(t)$ where $k$ is a finite extension of
$M(X,t)$ (respectively $k = M(X,t)$, if $t$ is a Galois
covering). Then, the field extension $K(X)/C(t)$ is defined over
$k$. By the usual dictionary between curves and function fields,
this means that Theorem (2.2) is proved.\\
For the proof of the above claim, we denote by $U(X,t,P)$ the
subgroup of $U(X,t)$ consisting of all $\sigma \in \Aut(C)$ such
that there is an isomorphism $f_\sigma: X^\sigma \ra X$ of curves
over $C$ such that the diagram
\[\xymatrix{
X^\sigma \ar[d]^{t^\sigma} \ar[rr]^{f_\sigma} && X \ar[d]^t \\
(\PP^1_C)^\sigma \ar[rr]^{\Proj(\sigma)} && \PP^1_C}\] commutes
and such that $f_\sigma(P^\sigma) = P$; here, $P^\sigma$ denotes
the point on $X^\sigma/C$ corresponding to $P$. Note that
$f_\sigma$ is unique since $\Aut(t)$ acts freely on the fibre
$t^{-1}(Q)$. Thus, mapping $\sigma$ to the automorphism of the
function field $K(X)$ induced by $f_\sigma$ yields an action of
$U(X,t,P)$ on $K(X)$ by $C$-semilinear field automorphisms which
fix $t \in K(X)$. Being the stabilizer of $[P]$ under the
(well-defined!) action $(\sigma, [P]) \mapsto
[f_\sigma(P^\sigma)]$ of $U(X,t)$ on $t^{-1}(Q)/\Aut(t)$, the
subgroup $U(X,t,P)$ has finite index in $U(X,t)$. If $t$ is a
Galois covering we in fact have $U(X,t,P) = U(X,t)$ since then
$t^{-1}(Q)/\Aut(t)$ has only one element. The meromorphic function
$z \in K(X)$ and hence the minimal polynomial of $z$ over $C(t)$
are invariant under the action of $U(X,t,P)$ defined above since
the image of $z$ under $\sigma \in U(X,t,P)$ has the same three
defining properties as $z$, as one easily checks. Now, Lemma~(1.6)
implies the above claim. Thus, the proof of Theorem~(2.2) is now
complete. \hspace*{\fill} $\Box$

{\bf (2.3) Remark.} Let $C=\bar{\QQ}$ and let $t$ be a Galois
covering. Then the second assertion of Theorem (2.2) can be proved
more quickly as follows.\\ There is obviously a model $t_L: X_L
\ra \PP^1_L$ of $t$ over a finite Galois extension $L$ of $\QQ$
such that $X_L$ has an $L$-rational point $P$ with $\QQ$-rational
and unramified image $Q:=t_L(P)$, such that all automorphisms of
$t$ are defined over $L$, and such that, for any $\sigma \in G:=
{\rm Image}(U(X,t) \ra \Aut(L))$, there is an isomorphism
$f_\sigma: X_L^\sigma \ra X_L$ of varieties over $L$ such that the
following diagram commutes:
\[\xymatrix{
X_L^\sigma \ar[d]^{t^\sigma_L} \ar[rr]^{f_\sigma} && X_L
\ar[d]^{t_L}
\\ (\PP^1_L)^\sigma \ar[rr]^{\Proj(\sigma)} && \PP^1_L.}\]
Since $\Aut(t_L)$ acts freely and transitively on $t_L^{-1}(Q)$,
there is a unique isomorphism $f_\sigma$ as above with
$f_\sigma(P^\sigma) = P$. Then we have $f_{\sigma \tau} = f_\sigma
\circ f_\tau^\sigma$ for all $\sigma \in G$. Now, the second
assertion of Theorem (2.2) follows from the comparatively
elementary Theorem (1.9).

\bigskip

\section*{3. The Theorem of Belyi}

We begin with the following proposition. It occurs implicitly at
several places in the literature and it is the analogue of a
well-known theorem in Algebraic Number Theory (e.g., see Theorem
(2.13) on p.~214 in \cite{Ne}). For the convenience of the reader
we include a short, self-contained and elementary proof which uses
only standard facts of the theory of Riemann surfaces and of the
theory of unramified topological coverings. Another proof using
triangle groups can be found in \S 1 in Wolfart's paper \cite{Wo}.

{\bf (3.1) Proposition.} {\em Let $S$ be a finite set of (closed)
points of the projective line $\PP^1_\CC$, and let $d\ge 1$ be a
natural number. Then there are at most finitely many isomorphism
classes of pairs $(X,t)$ where $X/\CC$ is a curve and $t:X \ra
\PP^1_\CC$ is a finite morphism of varieties over $\CC$ of degree
$d$ whose critical values lie in $S$.}

Here, two pairs $(X_1, t_1)$, $(X_2, t_2)$ as above are called
{\em isomorphic}, iff there is an isomorphism $f:X_1 \,\,
\tilde{\ra}\,\, X_2$ of varieties over $\CC$ with $t_2 \circ f =
t_1$.

{\em Proof.} By passing from a finite morphism $t: X\ra \PP^1_\CC$
to the continuous map $t(\CC): X(\CC) \ra \PP^1(\CC)$ between the
corresponding Riemann surfaces and by restricting $t(\CC)$ to the
preimage of the punctured sphere $\PP^1(\CC) \backslash S$, we
obtain a map from the set of isomorphism classes of pairs as above
to the set $\cM$ of homeomorphism classes of unramified
topological coverings of $\PP^1(\CC)\backslash S$ of degree $d$.
This map is injective. To see this, let $(X_1, t_1)$ and $(X_2,
t_2)$ be two pairs as above and let $g: X_1(\CC) \backslash
t_1^{-1}(S)\ra X_2(\CC) \backslash t_2^{-1}(S)$ be a homeomorphism
with $t_2(\CC) \circ g = t_1(\CC)$ on $X_1(\CC) \backslash
t_1^{-1}(S)$; then $g$ is biholomorphic, since
$t_i(\CC)|_{X_i(\CC) \backslash t_i^{-1}(S)}$, $i=1,2$, are
locally biholomorphic; by an elementary fact in Complex Analysis
(e.g., see Satz 8.5 on p.\ 48 in \cite{Fo}), the map $g$ can be
extended to a biholomorphic map $h: X_1(\CC) \ra X_2(\CC)$ with
$t_2(\CC) \circ h = t_1(\CC)$; now we apply the not very deep fact
that any biholomorphic map between complex curves is algebraic
(see section IV.11 in \cite{FK} or Lecture~9 in \cite{Be}) to get
an isomorphism $f: X_1 \ra X_2$ of varieties over $\CC$ with $t_2
\circ f  = t_1$; i.e, the pairs $(X_1, t_1)$ and $(X_2, t_2)$ are
isomorphic. Thus, it suffices to show that the set $\cM$ is
finite. Since any unramified topological covering of $\PP^1(\CC)
\backslash S$ is a quotient of the universal covering $p$ by a
subgroup of $\Aut(p) \cong \pi_1(\PP^1(\CC) \backslash S)$, we are
reduced to showing that there are at most finitely subgroups of
index $d$ of the fundamental group $\pi_1(\PP^1(\CC) \backslash
S)$. This follows from the facts that $\pi_1(\PP^1(\CC) \backslash
S)$ is finitely generated (in fact, $\pi_1(\PP^1(\CC) \backslash
S) \cong \langle \gamma_Q, Q\in S: \prod_{Q\in S} \gamma_Q =1
\rangle$ is a free group of rank $|S| -1$, see Aufgabe 5.7.A2 in
\cite{SZ}) and that a finitely generated group has only finitely
many subgroups of a given finite index (well-known and easy to
prove; it also follows from Theorem 7.2.9 on p.~105 in
\cite{Hal}). So, Proposition (3.1) is proved. \hspace*{\fill}
$\Box$

{\bf (3.2) Corollary}. {\em Let $X/\CC$ be a curve, let $t: X \ra
\PP^1_\CC$ be a finite morphism and let $K$ be a subfield of $\CC$
such that the critical values of $t$ are $K$-rational. Then the
moduli field of $t$ is contained in a finite extension of $K$.}

{\em Proof}. For any $\sigma \in \Aut(\CC/K)$, the critical values
of $t(\sigma): X^\sigma \,\, \stackrel{t^\sigma}{\longrightarrow}
\,\, (\PP^1_\CC)^\sigma \,\,
\stackrel{\Proj(\sigma)}{\longrightarrow}\, \, \PP^1_\CC$ lie in
$S$, too, and the degree of $t(\sigma)$ is the same as the degree
of $t$. So, by Proposition~(3.1), the orbit of the isomorphism
class of the pair $(X,t)$ under the obvious action of
$\Aut(\CC/K)$ is finite. Hence, the stabilizer is of finite index
in $\Aut(\CC/K)$. Furthermore, it is obviously contained in
$U(X,t)$. Now, Lemma~(1.4) and Lemma~(1.6) imply that the moduli
field $M(X,t) = \CC^{U(X,t)}$ is contained in a finite extension
of $\CC^{\Aut(\CC/K)} = K$. \hspace*{\fill} $\Box$

We are now ready to prove Belyi's theorem. The if-direction is a
consequence of Theorem (2.2) and Corollary (3.2) (see below). For
completeness sake, we also give a proof of the only-if-direction
(see Lemmas (3.4) through (3.6)).

{\bf (3.3) Theorem} (Belyi, 1979). {\em A complex curve $X$ is
defined over a number field, if and only if there exists a finite
morphism $t: X\ra \PP^1_\CC$ of varieties over $\CC$ with at most
$3$ critical values.}

{\em Proof.} First, we assume that there is a morphism $t: X \ra
\PP^1_\CC$ as above. After composing $t$ with an appropriate
fractional linear transformation, we may assume that the critical
values of $t$ lie in $S:=\{0, 1, \infty\}$. Then the moduli field
$M(X,t)$ is a number field by Corollary (3.2). Now, Theorem (2.2)
shows that $X$ is defined over a (may be, bigger) number field.
This proves the if-direction of Theorem (3.3). To prove the
only-if-direction, we introduce the notation $\Crit(f)$ for the
set of critical values of any morphism $f$ between curves. We
first choose an arbitrary morphism $t':X \ra \PP^1_\CC$ defined
over $\bar{\QQ}$ and apply Lemma~(3.4) below to $N:=\bar{\QQ}$ and
$t:= t'$, we then apply Lemma~(3.5) to $S:= \Crit(t')$ which
yields a certain morphism $p:\PP^1_\CC \ra \PP^1_\CC$, and we
finally apply Lemma~(3.6) to $T:= \Crit(p) \cup p(S)$ which yields
another morphism $q: \PP^1_\CC \ra \PP^1_\CC$. Then the
composition $t:= q \circ p \circ t'$ has at most 3 critical values
since, for any composition $g\circ f$ of morphisms between curves,
we obviously have $\Crit(g\circ f) = \Crit(g) \cup g(\Crit(f))$.
This completes the proof of Theorem~(3.3). \hspace*{\fill} $\Box$

{\bf (3.4) Lemma.} {\em Let $X/\CC$ be a curve defined over an
algebraically closed subfield $N$ of $\CC$ and let $t: X \ra
\PP^1_\CC$ be a finite morphism defined over $N$. Then the
critical values of $t$ are $N$-rational.}

{\em Proof.} Let $t_N: X_N \ra \PP^1_N$ denote a model of $t$ over
$N$, and let $\alpha: X \ra X_N$ and $\beta: \PP^1_\CC \ra
\PP^1_N$ denote the canonical projections. Then we have:
\begin{eqnarray*}
\lefteqn{{\rm Crit}(t) = t({\rm supp}(\Omega^1_{X/\PP^1_\CC})) =
t({\rm supp}(\alpha^*(\Omega^1_{X_N/\PP^1_N}))) } \\ && \subseteq
t(\alpha^{-1}({\rm supp}(\Omega^1_{X_N/\PP^1_N}))) \subseteq
\beta^{-1}(t_N({\rm supp}(\Omega^1_{X_N/\PP^1_N}))).
\end{eqnarray*}
This proves Lemma~(3.4) since the projection $\beta$ maps each
point of $\PP^1_\CC$ which is not $N$-rational to the generic
point of $\PP^1_N$. \hspace*{\fill} $\Box$

{\bf (3.5) Lemma.} {\em Let $S$ be a finite subset of $\bar{\QQ}$.
Then there is a non-constant polynomial $p \in \QQ[z]$ such that
$p(S)$ and the critical values of $p: \PP^1_\CC \ra \PP^1_\CC$ lie
in $\QQ \cup \{\infty\}$.}

{\em Proof.} We may and we will assume that $S$ is closed under
conjugation and use then induction on the number $n$ of elements
in $S$. If  $n\le 1$, we may take $p=z$. So, let $n > 1$. There is
a polynomial $p_1 \in \QQ[z]$ of degree $n$ such that $p_1(S) =0$
(namely the product of minimal polynomials of the elements in
$S$). We set $S_1:= p_1(\{r \in \bar{\QQ}: p'_1(r) = 0\})$. Then
$S_1 \cup \{\infty\}$ is the set of critical values of $p_1$,
$S_1$ has at most $n-1$ elements, and $S_1$ is closed under
conjugation again. By the induction hypothesis, there is a
polynomial $p_2 \in \QQ[z]$ such that $p_2(S_1)$ and the critical
values of $p_2$ lie in $\QQ \cup \{\infty\}$. Then the
composition $p:=p_2 \circ p_1$ satisfies
\[\Crit(p) = \Crit(p_2)
\cup p_2(\Crit(p_1)) = \Crit(p_2) \cup p_2(S_1 \cup \{\infty\})
\subseteq \QQ \cup \{\infty\}\] and $p(S) = p_2(p_1(S)) =
p_2(\{0\}) \subseteq \QQ$, as desired. \hspace*{\fill} $\Box$

{\bf (3.6) Lemma.} {\em Let $T$ be a finite subset of $\QQ$. Then
there is a non-constant morphism $q: \PP^1_\CC \ra \PP^1_\CC$ such
that $q(T)$ and the critical values of $q$ lie in
$\{0,1,\infty\}$.}

{\em Proof.} We use induction on the number $r$ of elements in
$T$. If $r \le 3$, there is  a fractional linear transformation
$q$ with $q(T) \subseteq \{0, 1, \infty\}$. So, let $r > 3$. After
composing with an appropriate fractional linear transformation, we
may assume that $0, 1, \infty \in T$ and that there is a fourth
point in $T$ which lies in the interval between $0$ and $1$, i.e.,
which is of the form $\frac{m}{m+n}$ where $m,n \in \NN$. We now
consider the polynomial
\[q_1(z) := \frac{(m+n)^{m+n}}{m^m n^n}
z^m (1-z)^n \in \QQ[z].\] Then $q_1$ maps the set of four points
$0, \frac{m}{m+n}, 1, \infty$ onto $\{0,1,\infty\}$ and the
critical values of $q_1$ lie in $\{0,1, \infty\}$ since the
derivative $q'_1(z)$ equals $z^{m-1} (1-z)^{n-1}((m+n)z-m)$ up to
a constant. By the induction hypothesis applied to $q_1(T)$, there
is a morphism $q_2: \PP^1_\CC \ra \PP^1_\CC$ such that
$q_2(q_1(T))$ and the critical values of $q_2$ lie in
$\{0,1,\infty\}$. Then the composition $q := q_2 \circ q_1$
satisfies
\[\Crit(q) = \Crit(q_2) \cup q_2(\Crit(q_1)) = \Crit
(q_2) \cup q_2(\{0,1,\infty\}) \subseteq \{0,1,\infty\}\] and
$q(T) = q_2(q_1(T)) \subseteq \{0,1,\infty\}$, as desired.
\hspace*{\fill} $\Box$

The following statement is a corollary of our proof of the
if-direction of Belyi's theorem. It gives a bound for the degree
of the field of definition of a complex curve which allows a
finite morphism to the projective line with at most 3 critical
values. For this, let $M_d$ denote the number of subgroups of
index $d$ in a free group of rank 2. We have the following
recursion formula for $M_d$:
\[M_d = d (d!) - \sum_{i=1}^{d-1} (d-i)! M_i\]
(see Theorem 7.2.9 on p.\ 105 in \cite{Hal}).

{\bf (3.7) Corollary.} {\em Let $X/\CC$ be a curve and let $t: X
\ra \PP^1_\CC$ be a finite morphism of degree $d$ with ${\rm
Crit}(t) \subseteq \{0,1, \infty\}$. Let $a$ denote the number of
elements in $\Aut(t)$. Then $X$ and $t$ are defined over a number
field $K$ with $[K:\QQ] \le \frac{d}{a}M_d$.}

{\bf Proof.} The fundamental group of the punctured sphere
$\PP^1(\CC) \backslash \{0, 1, \infty\}$ is a free group of rank~2
(see the proof of Proposition (3.1)). Thus, the orbit of the
isomorphism class of the pair $(X,t)$ under the action of
$\Aut(\CC)$ has at most $M_d$ elements (see the proof of
Proposition (3.1)). This implies that $[M(X,t) :\QQ] \le M_d$ as
in the proof of Corollary (3.2). By the proof of Theorem (2.2),
$X$ and $t$ are defined over $\CC^{U(X,t,P)}$ and the index of the
subgroup $U(X,t,P)$ in $U(X,t)$ is less than or equal to
$\frac{d}{a}$. Now, Lemma~(1.6) proves Corollary (3.7); note that
$U(X,t)$ is closed by Lemma~(1.5) and Theorem (2.2).
\hspace*{\fill} $\Box$

{\bf (3.8) Remark}. In Corollary (3.7), the number $M_d$ may of
course be replaced be the (smaller) number of isomorphism classes
of pairs $(X',t')$ where $X'$ is a curve and $t': X' \ra
\PP^1_\CC$ is a finite morphism of degree $d$ such that
$\textrm{Crit}(t') \subseteq \{0,1,\infty\}$ and such that, in
addition, the ramification indices of $t'$ are the same as those
of $t$. Besides the ramification indices, one may also take into
account further Galois invariants of a Belyi surface (like the
monodromy group or the cartographic group). It would be
interesting to get explicit formulas for these (much) sharper
bounds or to get at least explicit estimations which substantially
improve the bound $M_d$.

\bigskip

\bigskip

Faculty of Mathematical Studies\\ University of Southampton\\
Southampton SO17 1BJ\\ United Kingdom\\ {\em e-mail:}
bk@maths.soton.ac.uk.

\end{document}